\input amstex.tex
\input amsppt.sty

\def\a{\alpha}

\def\g{\gamma}
\def\d{\delta}
\def\s{\sigma}
\def\V{\Cal V}

\magnification=1200
\tolerance=5000
\refstyle{A}

\topmatter

\abstract The infinite-dimensional Clifford algebra has a maze of inequivalent irreducible unitary representations. Here we determine their type -real, complex or quaternionic. Some, related to the Fermi-Fock representations,  do not admit 
any real or quaternionic structures. But there are many 
 on $L^2$ of the circle that do and which seem to have analytic meaning.
\endabstract

\author 
{Esther Galina$^{1,3}$, Aroldo Kaplan$^{1,2,3}$, 
Linda Saal$^{3}$}
\endauthor

\affil 
$(^1)$ CIEM-CONICET, Argentina,  $(^2)$  Department of Mathematics and Statistics, University of Massachusetts, Amherst, $(^3)$ FaMAF, Universidad Nacional de C\'ordoba
\endaffil

\leftheadtext{E. Galina, A.Kaplan, L. Saal}

\title Reality of non-Fock Spinors \endtitle

\thanks
This research was supported in part by CONICET, FONCYT, CONICOR and UNC.
\endthanks
\endtopmatter
\document
\subhead{Table of contents}
\endsubhead 

1. Introduction

2. Spinors as dyadic objects

3. Real and quaternionic structures 

4. $L^2(\Bbb T)$ as a spinor space 

5. Dyadic difference operators

6. Kaplansky's division algebras 

\

\subhead 1. Introduction
\endsubhead

Let $H$ be a separable real prehilbert space and  let $C(H)$ be the Clifford algebra of $H$, i.e.,  the quotient of the tensor algebra $\Cal T(H)$ of $H$ by the ideal generated by the elements of the form
 $$h \otimes h'+h'\otimes h + 2<h ,h'>\qquad   h,h'\in H.$$ 
Here we parametrize all the equivalence classes of representations of $C(H)$ on a separable real Hilbert space $U$, where 
$H\subset C(H)$ acts via skew-symmetric operators (``orthogonal"). 
$U$ is a {\it space of (real) spinors}. In infinite dimensions there 
is ``a true maze" of inequivalent irreducible ones, in striking contrast 
to the finite case.

Choosing an orthonormal basis of $H$ and letting $J_k$ denote the action of $k^{\roman th}$ element of the basis on $U$, the operators
 $J_1, J_2, ... ,$ are orthogonal complex 
structures on $U$ which anticommute with each other. 
We will often ignore $H$ altogether and regard a spinor structure 
on a real Hilbert space $U$ as 
a sequence of linear operators 
 $J_1, ..., J_k, ....,$
on it satisfying
$$ ||J_kv||=||v||, 
\qquad J_k^2=-I, \qquad 
J_kJ_l=-J_lJ_k$$
for all $k$ and all $l\not=k$.

The  complex spinors, i.e., the unitary representations of $\Bbb C\otimes C(H)$, are the same as the representations of $C(H)$ on separable complex Hilbert spaces $V$ satisfying
$$||h\cdot v|| = ||h||\ ||v||$$
for $h\in H\subset C(H)$ and $v\in V$. As we explain in \S 2, when $\dim H$ is even or infinite, $\Bbb C\otimes C(H)$ is the same as  $\Bbb C\otimes \Cal T(H)$ modulo the so-called Canonical Anticommutations Relations of Quantum Field Theory. G\aa rding and Wightman's 
parametrized the representations of the latter in [GW1]. Therefore, their 
parametrization yields a corresponding one 
of all the complex spinor structures up to unitary equivalence. We describe the result in detail in \S 2.

To parametrize the real representations it is then  
enough to determine the values of the G\aa rding-Wightman (or GW) parameters whose 
corresponding complex representation admits a  invariant real  structure, i.e., 
a $\Bbb C$-antilinear, norm-preserving operator $S$ 
such that 
$$S^2=I.$$
Then 
$$U=\{v: \ Sv=v\} $$
is an invariant real form of $V$ which, by restriction, provides a real representation 
of $C(H)$. 
Since the complexification of a real representation of 
$C(H)$ is a representation of $\Bbb C\otimes C(H)$, 
one obtains a parametrization of the former.

Similar arguments yield those complex representations posessing a $\Bbb C$-antilinear, norm-preserving operator $Q$ such that 
$$Q^2=-I.$$

In the classical (finite dimensional) case, a complex representation of $C(H)$ is defined as of {\it real, quaternionic} or {\it complex  type}, according to 
whether it admits an $S$, a $Q$, or neither, 
conditions that are mutually exclusive when the representation 
is irreducible.  In the physics literature
$S$ and $Q$ are  called  
{\it charge conjugation  operators} and the irreducible representations of real type 
{\it Majorana spinors}. 
By the way, 
we recover the classical result of Cartan and Dirac, namely that the unique irreducible complex representation of $C(\Bbb R^{2m})$  is of real type 
if and only if $m\equiv 0,3$ modulo $4$ and of quaternionic type otherwise.  
The geometric and physical significance of this in the finite dimensional case is well 
known (see, e.g., [C] and other chapters in the same book).

In infinite dimensions we
find {\it mazes of inequivalent  irreducible spinors  
of each of the three types}. Those corresponding to the Fermi-Fock representations of
the Canonical Anticommutation Relations admit no real or quaternionic structures. But there are plenty of natural irreducible 
spinor structures of real and of quaternionic type on $L^2$ of the circle. 

It is important to note that the final result is 
just  a {\it parametrization} of {\it all} 
equivalence classes of unitary spinor structures of the 
three types.  The questions of irreducibility and equivalence of the GW representations for the various values of the parameters are not completely resolved yet. 
Also, although 
every unitary representation of this algebra 
is completely reducible, the reduction is highly non-unique. These constitute insurmountable obstacles for proving most general statements about spinors in infinite dimensions using the GW parametrization. In a way, one purpose of this paper is to
exhibit one problem, namely the classification into types, for which these obstacles can be surmounted and has a neat answer.

As to other purposes, we mention some preliminary algebraic and analytic consequences. 

First, and much like when
$\dim H =1,3,7$, such real representations are 
in correspondence with certain 
 {\it division algebras} 
considered by Kaplansky [K] 
as hypothetical infinite-dimensional analogs of the
Octonions -although he himself was doubtful of their existence. Here is then 
a complete parametrization of such algebras. 
Of course, they are not commutative or associative, 
and they contain only one-sided inverses. but there are mazes of inequivalent ones. Regardless of their numerical status, {\it their 
automorphism groups are reductive and come unitarily represented}. By restriction, 
one obtains mazes of new irreducible unitary representations of  
the  classical infinite dimensional Lie groups and algebras that appear as factors.

Secondly, there are interesting families of representations of $C(H)$ on $L^2(\Bbb
T)$ (or
$L^2(\Bbb R)$) of real or quaternionic type
which have some analytic  content.  For example, those of real type yield all 
manners of fitting the  Hilbert transform 
$$\Cal Hf(x) =\int_{-\infty}^{\infty} 
{f(y)\over x-y}\ dy$$
into a sequence of {\it mutually anticommuting real singular integral operators $\Cal H, \Cal H_2, \Cal H_3, ...$
of square $-I$}. The corresponding kernels 
are dyadic twistings of the Hilbert kernel and lead to analogs of the 
Cauchy kernel. In the quaternionic case, 
the charge conjugation operators are obtained by twisting real ones by Haar's mother wavelet. Intriguing as they may be, we will not go into much detail about these issues here. 

Instead, we will discuss  two operators, 
$$D = \sum_{k=1}^{\infty} a_k
\partial_k, \qquad D' = \sum_{k=1}^{\infty} a_k^*
\partial_k$$
where $a_k$, $a_k^*$, are the creation and anihilation 
operators associated to any spinor structure and 
the $\partial_k$ are certain dyadic difference operators. Notably, for the standard  Fermi-Fock representations 
they diverge off the vacuum. But for the  spinor structures in
$L^2(\Bbb T)$ that we dicuss below, they have a dense domain  and are mutually conjugate
under any charge conjugation operator.  We found 
remarkable that for one of 
these families, parametrized by infinite  matrices of $0$'s and $1$'s, the associated operators
$D$ and
$D'$, which are far from self adjoint, {\it can be diagonalized over
$\Bbb Z$:  with integral eigenvalues and   eigenfunctions that
 are
 polynomials with integral coefficients} in the classical periodic Rademacher functions.

Any connection of all this with the real world 
must take into account that,
as we prove below, a
real  or quaternionic structure requires that, 
in the standard statistical intepretation, changing all the occupied states to
non-occupied and viceversa be a well defined 
operation. This may be an unlikely feature  for 
fermions, but not necessarely for other systems of $0$'s and $1$'s. Indeed, the properties of 
$D,D'$ and the higher Hilbert transforms $\Cal H_k$
seem more related to wave packings, splines and binary codes than to any 
particles or fields. 

As this is a preliminary version, some proofs 
are only sketched, others are found in [GKL] and 
a complete version will be ready shortly.

We thank H. Araki, J. Baez, A. Jaffe, and A. Kirillov for their useful advise.

\

\subhead 2. Spinors as dyadic objects
\endsubhead

Let
$$X=\Bbb Z_2^{\infty}$$
be the set of sequences $x=(x_1,x_2,\dots)$ of $0$'s and $1$'s, and  $\Delta\subset X$ the
subset consisting of sequences with only finitely many $1$'s. Then $X$ is an  abelian
group under componentwise addition modulo $2$ and $\Delta$ is the subgroup generated by
the sequences $\delta^k$, where $\delta^k_j$ is the Kronecker symbol.  
The product topology  
on $X$ is compact and is generated
by the sets 
$$X_k=\{x\: x_k=1\},\qquad 
X_k'=\{x\: x_k=0\},$$ 
which, therefore, also generate the canonical $\s$-algebra 
of Borel sets in $X$.

We will realize all the complex spinor 
structures  on 
$L^2$ spaces of $\Bbb C$-valued functions on $X$ or 
direct integrals thereof. As a  motivation,  
let us realize the standard finite  even-dimensional spinors in this manner.  For each positive integer $N$ consider the vector space 
$$V_N= \Bbb C^{\ \Bbb Z_2^N}.$$ 
Then, clearly, $\dim V_N = 2^N$ and the operators
$$
\aligned
J_k f(x)&= -i(-1)^{x_1+\dots +x_{k-1}} \ f(x+\delta^k)\\
J_k' f(x)&=  \ (-1)^{x_1+\dots  +x_k} \ f(x+\delta^k)\\
\endaligned
\tag{2.1}
$$
where $1\leq k\leq N$, $x\in \Bbb Z_2^N$, addition 
is modulo 2 and the $\d^k$ is the standard 
basis of $\Bbb Z_2^N$, define an irreducible  complex representation of 
the Clifford algebra $C(\Bbb R^{2N})$ -the unique 
one modulo equivalence. In spite of its simplicity 
and of being implicit 
in the work of Friedrichs, G\aa rding, Wightman and 
von Neumann on the Anticommutation Relations, this description of even dimensional spinors 
does not seem to have been made explicit or exploited before.

The unitarity underlying the finite case is relative 
to the natural $L^2$ inner product in $V_N$, which in turn is associated to the measure on  
$\Bbb Z_2^N$ where each point has measure $1$. 
But changing the measure to any equivalent one 
or changing the target space of the functions does not 
change the equivalence class of the 
representation.

This is no longer so when $N=\infty$: in order to reach 
all equivalence classes  one 
must allow for more general measures on the group 
$X=\Bbb Z_2^\infty$ and replace $\Bbb C$-valued functions for sections of appropriate 
fiber spaces over $X$. 
Two canonical, but very different measures on $X$ that generalize 
the  finite case are:

\smallskip

- $\mu_X$, the Haar measure of $X$.

- $\mu_\Delta$, concentrated in the discrete set $\Delta$,
where $\mu_{\Delta}(\{\delta\} )=1$. $\mu_\Delta$ could be called the Fermi-Fock 
measure.
\smallskip
The first is invariant under  all translations in $X$ 
while the second is invariant only under those from
$\Delta$. It is $\mu_\Delta$ 
that leads to the representations that appear most in QFT, however implicitely. It 
 ignores all the points $x$ with infinitely many 
$x_i=1$, or occupied states, on the basis that the total number of particles -fermions in this case, must be finite. In any case, (2.1) define irreducible representations of $C(H)$ on $L^2(X,\mu_X)$ and 
$L^2(X,\mu_\Delta)$, respectively, which we will prove to be inequivalent.

\

Recall that two  measures $\lambda, \mu$ on a the same Borel 
algebra of sets can be said to be {\it equivalent} 
if there exists locally integrable functions, denoted by $d\lambda/d\mu$ and 
$d\mu/d\lambda$, such that for any measurable $A$, 
these Radon-Nykodim derivatives satisfy
$$\lambda(A)=\int_A {d\lambda\over d\mu} \ d\mu,
\qquad \mu(A)=\int_A {d\mu\over d\lambda} \ d\lambda.
$$ 
$\mu$ is said to be {\it quasi-invariant} under a set of transformations $\{T\}$ of the 
underlying space, if the translated measures 
$$\mu_T(A):=\mu(T(A))$$
are all equivalent to $\mu$. 
Now, consider triples 
$$(\mu,\V,\Cal C) $$ 
where

\

{\it $\bullet\ \ \mu$ is a positive  Borel measure on $X$, 
quasi-invariant under $\Delta$.

\

$\bullet\ \V= \{V_x\}_{x\in X}$ is a family of complex 
Hilbert spaces a.e. invariant under 
translations by $\Delta$ and 
such that the function
$$\nu(x)=\dim V_x$$
is measurable.

\

$\bullet\ \ \Cal C = \{c_k(x):\ k\in \Bbb Z_+,\ x\in X\}$   is a family of unitary operators 
$$c_k(x): V_x \rightarrow V_{x+\delta^k}= V_x$$
depending measurably on $x$ and satisfying 
$$
\aligned
c^*_k(x)&=c_k(x+\delta^k)\\
c_k(x)c_l(x+\delta^k)&= c_l(x)c_k(x+\delta^l)
\endaligned
\tag {2.2}
$$
for all $\delta\in \Delta$ and almost all $x\in X$. }

\

We will often write 
 $(\mu,\nu,\Cal C) $ instead of 
 $(\mu,\V,\Cal C))$, in view of the fact that 
changing $\V$ unitarily will yield equivalent 
representations.
Given such triple, consider the Hilbert 
space 
$$V=V(\mu,\nu,\Cal C)=\int_X^\oplus V_x \ d\mu(x).
$$
For example, when $V_x=\Bbb C$ for all $x$,
$$V(\mu,1,\Cal C)= L^2(X,\mu),$$
the ordinary $L^2$ space of $\Bbb C$-valued functions. 
Finally, define operators on $V$ by
$$
\aligned
J_k f(x)&= -i(-1)^{x_1+\dots +x_{k-1} }\  c_k(x) \ \Big( \sqrt{\frac {d\mu(x+\delta^k)}{d\mu(x)}} \ f(x+\delta^k)\Big)\\
J_k' f(x)&= \ \ (-1)^{x_1+\dots +x_{k}} \ \ c_k(x) \ \Big(\sqrt{\frac {d\mu(x+\delta^k)}{d\mu(x)}} \ f(x+\delta^k)\Big)\\
\endaligned
\tag {2.3}
$$
where an $f\in V$ is regarded as an assignement $x\mapsto f(x)\in V_x$ and all $+$ are modulo 2. 

\proclaim{Theorem 2.4} The operators $J_1,J_1',J_2,J_2',...$ are mutually anticommuting orthogonal complex structures and, therefore, 
define a (complex) spinor structure on $V$. Conversely, 
every spinor structure on a separable Hilbert space is unitarily equivalent to some $V(\mu,\nu,\Cal C)$.
\endproclaim

For the proof, one  observes that the 
Clifford commutation relations for the $J$'s translate into the Fermi commutation relations for the operators
$$a_k={1\over 2}\big(iJ_k-J'_k\big),$$
 which, according to [GW1], are themselves
parametrized up to equivalence by the triples $(\mu,\nu,\Cal C)$.
The Fermi-Fock representation corresponds to the triple  $(\mu_{\Delta},1,\{1\})$.              
Von Neumann's first examples of non-Fock representations, were 
infinite tensor products, which in our notation 
are the   $$V(\mu_X,1,\Cal C^\otimes)$$
 where
$\mu_X$ is the Haar measure on $X$ and 
$$c_k^{\otimes}(x) =\omega_k^{(-1)^{x_k}},$$
the $\omega_k$ being arbitrary complex numbers of 
absolute value $1$. In particular,  
$$V(\mu_X,1,\{1\}).$$
 is one such.

\subhead 
3. Real and Quaternionic structures
\endsubhead

If $U$ is a real module over $C(H)$, then $\Bbb C\otimes U$ is a complex module over $\Bbb C\otimes C(H)$, which comes with the 
$C(H)$-invariant decomposition
$$\Bbb C\otimes U = U \oplus_{\Bbb R} iU.$$
$U$ is an {\it invariant real form} of $\Bbb C\otimes U$. Conversely, any complex module over $C(H)$ with an invariant real form determines a real module over $C(H)$. Hence, by determining all the invariant real 
forms of the G\aa rding-Wightman modules we will be 
parametrizing all the real representations of $C(H)$ 
up to orthogonal equivalence.

The first problem is equivalent to determining all the charge-conjugation operators of the 
representations 
$V(\mu,\nu,\Cal C)$, i.e., the $\Bbb C$-antilinear operators $S:V\to V$ which
commute with the action of $C(H)$ and such that
$$
S^2=1, \qquad ||Sf|| = ||f||. \tag {3.1}
$$
The invariant real form in question is then $\{v\in V:\ Sv=v\}$.

Let 
$$x\mapsto \check x$$
be the involution of $X$ which 
changes all $0$'s to $1$'s and viceversa. Modulo 2,
$$\check x= x+ \bold1$$
where $\bold1_k=1$\ $\forall k$. We have induced involutions on subsets of $X$ and on functions and 
measures on $X$:
$$\check A=\{\check x: \ x\in A\}, \qquad 
\check f(x) = f(\check x), \qquad 
\check \mu(A) = \mu(\check A).$$

\proclaim{Theorem 3.2} $V(\mu,\nu,\Cal C)$ admits an invariant real form if and only if 
the measures $\mu$ and $\check\mu$ are equivalent, $\check\nu(x) =\nu(x)$ for almost all
$x\in X$  and there exist a measurable family of operators 
$$r(x)\:V_x \to V_{\check x}\cong V_x$$
which are $\Bbb C$-antilinear, preserve the norm and 
satisfy
$$
\eqalign{r(x)r(\check x)&=1\cr
 r(x)c_k(\check x)&=(-1)^{k}c_k(x)r(x+\delta^k)\cr}
\tag {3.3}$$
for all $k\in \Bbb N$ and almost all $x\in X$.
\endproclaim

Sketch of proof: let
$$
(Tf)(x)= \sqrt{\frac{d\tilde\mu(x)}{d\mu(x)}} f(\check x).
$$
If $S$ is  an invariant real structure the the product $TS$ must commute
with the operators $a_ka_k^*$ and $a_k^*a_k$ for all $k$. This is a commuting set of
projections  for which 
$V=\int_X^\oplus V_x \ d\mu(x)$ is the spectral decomposition. From this one can deduce
that
$TS$ must act pointwise on each $V_x$.  This is our $r(x)$.

\

All this applies to the finite case as well. If $\dim H=2m$,  $\mu$ and $\check\mu$ are
equivalent for any $\mu$. From (3.3) one deduces 
$$r(\bold1)=(-1)^{\frac {m(m+1)}2}r(\bold 0).$$ 
Assuming, as we may, that $r(\bold 0)$ is the standard conjugation on $\Bbb C$, we see that $V$ splits 
over $\Bbb R$ if and only if  $m(m+1)/2$ is an even integer, i.e., for 
$$m\equiv 0,3\ \ (\roman {mod}\ 4)$$
as is well known.

Assume now that $V$ is infinite dimensional and separable. The axiom of choice implies that 
there are always plenty of solutions $r(x)$ 
to the equations (3.3) but 
most of them -and often all,  are non-measurable. Indeed, the latter turns out to be 
the case in the following two cases.
 
\proclaim{Corollary 3.4} If  $\mu$ is discrete and $V$ is irreducible over $\Bbb C$, then it is irreducible over $\Bbb R$. In particular, this is the case for the Fermi-Fock representations.
\endproclaim

\proclaim{Corollary 3.5} The tensor product representations
$V(\mu_X,1,\Cal C^{\otimes} )$ are  irreducible over $\Bbb R$. 
\endproclaim

The proofs of these results involve arguments of ergodicity.

\

Next, we describe a {\it standard} form for representations 
of real type for which the operator-valued function 
$r(x)$ is constant, namely
$$r(x)v = \bar v$$
with respect to a fixed choice of real form 
in  each integrand $V_x$. A maze 
of examples will come out of it. 

Let 
$$U=\int_X^{\oplus}U_x \ d\mu(x)$$ be a direct integral of {\it real}  Hilbert spaces satisfying 
$$U_{x+\delta} = U_x,\qquad U_{\check x} = U_x$$
for all $\delta\in \Delta$ and almost all $x\in X$. With $V = {\Bbb C}\otimes U$ and $V_x={\Bbb C}\otimes U_x$,
$$V = \int_X^{\oplus}V_x \ d\mu(x).$$ 
Clearly,  $V_{x+\delta} = V_x=V_{\check x}$, $U_x$ is a  real form of $V_x$ and $U$ one of $V$. Denote by $\ \bar{}\ $ the corresponding conjugations. If $f\in V$ and 
$A\in End(V)$ set  
$$\bar f (x):= \overline{f(x)},\qquad 
\bar A (f) := \overline{A(\bar f)}.\leqno{(3.6)}$$ 
$U$ is not the invariant real form we are looking for -this would be incompatible with 
the anticommutation relations. Instead we have 

\proclaim{Theorem 3.7} If
 $\check\mu$ is equivalent to $\mu$, $\check \nu =\nu $  and the operators $c_k$ satisfy 
$$\overline{c_k(x)} = (-1)^k c_k(\check x),$$ then
$$V^{\Bbb R} = \{f\in V: \overline{f(x)} = \sqrt{\frac{d\check\mu(x)}{d\mu(x)}} {f(\check x)}\}$$
is an invariant real form of 
$V=V(\mu,\nu,\Cal C).$ 
\endproclaim

\proclaim{Theorem 3.8} 
Any unitary representation of $C(H)$ with an invariant real structure is unitarely equivalent to one in standard form.
\endproclaim

Perhaps the simplest  infinite-dimensional Majorana 
spinors are given by  $V(\mu_X,1,\Cal C)$ with 
$\mu_X$ being the Haar measure of $X$ and
the $c_k$ given by the dyadic Rademacher functions 
$$\eqalign{c_{2\ell}(x)&=1\cr
c_{4 \ell+1}(x) &= (-1)^{x_{4\ell+3}}\cr
c_{4\ell+3}(x) &= (-1)^{x_{4\ell+1}}.\cr}$$

\proclaim{Theorem 3.9} With these $c_k$, $V(\mu_{X},1,\Cal C)$ is irreducible 
over $\Bbb C$, but  the real form 
$$
L^2(X)^{\Bbb R} = \{f\in L^2(X): \ f(\check x)= \overline{f(x)} \}
$$
is an invariant real subspace. The  real representation 
on $L^2(X)^{\Bbb R}$ so obtained, is irreducible over
$\Bbb R$ and  does not arise from any representation of $\Bbb C\otimes C(H)$ by
restriction of the scalars. 
\endproclaim

\

The quaternionic case is treated similarly, although in the irreducible case   the  operator-valued function $q(x)$ 
cannot be taken to be constant and does not arise 
from any quaternionic structure in each $V_x$, like 
real ones do. It is not obvious 
{\it a priori} that $V(\mu,\nu,\Cal C)$ can support any 
quaternionic structure when 
$\nu(x)=1$.

\proclaim{Theorem 3.10} $V(\mu,\nu,\Cal C)$ admits an invariant quaternionic 
structure if and only if $\mu$ and $\check\mu$ are equivalent, $\check \nu(x) =\nu(x) $ for almost all $x\in X$  and there exist a measurable family of operators 
$$q(x)\:V_x \to V_{\check x}\cong V_x$$
which are $\Bbb C$-antilinear, preserve the norm and 
satisfy
$$
\eqalign{q(x)q(\check x)&=-1,\cr 
q(x)c_k(\check x)&=(-1)^{k}c_k(x)q(x+\delta^k)\cr}\tag
{3.6}$$ for all $k\in \Bbb N$ and almost all $x\in X$.
\endproclaim

\proclaim{Corollary 3.11} If  $\mu$ is discrete and $V$ is irreducible over $\Bbb C$, then $V(\mu,\nu,\Cal C)$ admits no quaternionic structure. In particular, the Fermi-Fock representations are of complex type.
\endproclaim

\proclaim{Corollary 3.12} The tensor product representations
$V(\mu_X,1,\Cal C^{\otimes})$ 
 are all of complex type. In particular, 
this is so for  $V(\mu_X,1,\{1\})$.
\endproclaim

Finally, we give a standard form for spinors of 
quaternionic type, for which the operator-valued 
function $q$ will be
$$q(x)v = (-1)^{x_1}\bar v.$$
The bar indicates that we are in the context 
of (3.6), where $V$ comes with a (non-invariant) 
real form $U$. With this understood we have

\proclaim{Theorem 3.13} If
 $\check\mu$ is equivalent to $\mu$, $\check\nu(x) =\nu (x)$ and the operators $c_k$ satisfy 
$$\overline{c_1(x)} = c_1(\check x),$$
$$\overline{c_k(x)} = (-1)^{k} c_k(\check x),$$
$\forall k\geq 2$ and almost all $x\in X$, then
$$Qf(x)= (-1)^{x_1} \sqrt{\frac{d\mu(\check x)}{d\mu(x)}} 
\overline{f(\check x)}$$
is an invariant quaternionic structure in $V(\mu,\nu,\Cal C).$ 
\endproclaim

\proclaim{Theorem 3.14} Any
unitary representation of $C(H)$ with an invariant quaternionic structure  is unitarily equivalent to 
a standard one.
\endproclaim

The simplest  irreducible infinite-dimensional  
spinors of quaternionic type  are 
realized in $L^2(X,\mu_X)$ as $V(\mu_X,1,\Cal C)$, 
with 
$$\eqalign{c_{2\ell}(x)&=1=c_1(x)\cr
c_{2 \ell+1}(x) &= (-1)^{x_{2\ell+3}}\ \ (\ell\geq 1)\cr
c_{2\ell+3}(x) &= (-1)^{x_{4\ell+1}}
\ \ (\ell\geq 1).\cr}$$

The following {\it dyadic representations}
give many examples of real and quaternionic spinors 
with special properties. Recall that
the {\it Walsh  functions}, 
as functions on $X$, are defined by
$$\phi_\alpha(x)=(-1)^{\sum \alpha_k x_k}$$
for $\alpha\in\Delta$, which are precisely the characters 
of $X$. Setting
$$\s^k = \d^1+\cdots + \d^k$$
then
$$\phi_{\s^k}(x)=(-1)^{x_1+\cdots +x_k}$$
which appear as  multipliers in the definition of 
the operators $J,J'$.

Let $\Gamma$ denote the set of infinite symmetric matrices $\g$ of $0$'s and $1$'s, with 
 only finitely 
many $1$'s in each row  or column and 
none along the diagonal. We regard each row $\g^k$ as a point in $\Delta$. $\Gamma$ contains the disjoint 
subsets
$$
\Gamma_{1} = \{\g\in\Gamma: \ \
 \sum_j \g^k_j\equiv k\  \ \ \forall k\}$$
$$
\Gamma_{-1} = \{\g\in\Gamma: \ \ 
\sum_j \g^1_j\equiv 0,\ 
\sum_j \g^k_j\equiv k\ \ \ \forall k\geq 2\}
$$
where the congruences are modulo 2. In other words, $\Gamma_1$ consists of the matrices 
where the  number of $1$'s in a row 
has the same  parity as the position of that 
row, while for $\Gamma_{-1}$ the condition is 
the same except
for the first row, for which it is reversed. 

In what follows we will take $\mu=\mu_X$, the Haar measure on $X$ and $\nu=1$, so that
$$V=L^2(X,\Bbb C).$$

\proclaim{Theorem 3.15} {\it For any $\g\in \Gamma$, the multiplier operators 
$$c_k(x) = \phi_{\g^k}(x)$$
 satisfy (2.2) and, therefore, 
$$J_kf(x) = -i\phi_{\s^{k-1}+\g^k}(x)f(x+\d^k)$$
$$J_k'f(x) = \phi_{\s^{k}+\g^k}(x)f(x+\d^k)$$
define a spinor representation. This 
representation is irreducible /$\Bbb C$. 
If $\g\in \Gamma_1$ (respectively, 
$\g\in \Gamma_{-1})$ then the representation is 
of real 
(respectively, quaternionic) type, in its respective 
standard form. Those of real type have 
$$L^2(X)^{\Bbb R} = \{f\in L^2(X): \ \overline{f(x)}=f(\check x)\}$$
as the invariant real form.
\endproclaim

On the Walsh 
basis, 
$$\eqalign{J_k\phi_\a &= -i(-1)^{\a_k} \phi_{\a+ \gamma^k +\sigma^{k-1}}\cr
 J_k'\phi_\a &= (-1)^{\a_k} \phi_{\a+ \gamma^k +\sigma^{k}}\cr}$$

\subhead 
4. $L^2(\Bbb T)$ as a space of spinors 
\endsubhead 

Any representation $V(\mu_X,1,\Cal C)$ where $\mu_X$ 
is the Haar measure and $\nu=1$,
can be realized on the standard $L^2(\Bbb T)$ of complex-valued functions on the circle. 

Indeed, we 
can identify each $V_x$ with $\Bbb C$, so that 
$V= L^2(X)$ and now use the dyadic expansions to identify $X$ with the interval $(0,1)$, or with $\Bbb T$ -except for a set of measure zero. The Haar 
measure on $X$ becomes the Lebesgue measure on $(0,1)$, and the Haar measure on $\Bbb T$. This is in spite of the fact that translations in  $X$ do not correspond to rigid rotations of $\Bbb T$ and that, as a topological space, $X$ is 
homeomorphic, not to $\Bbb T$, but to the Cantor set. 
The operation $x\mapsto \check x$ in $X$ correspond to $y\mapsto 1-y$ in $(0,1)$ which, on $\Bbb T\subset \Bbb C$, becomes ordinary complex conjugation.  
In particular,  for a standard representations of real type with $\nu=1$,  the 
$c_k$'s are functions from $\Bbb T$ to itself satisfying
$$\overline{c_k(t)} = (-1)^k c_k(\bar t)$$
and the $C(H)$-invariant real form is
$$
L^2(\Bbb T)^{\Bbb R} = \{f\in L^2(\Bbb T): \overline{f(t)} = {f(\bar t)}\}
$$

Via $X\approx\Bbb T$ the  functions 
$\phi_\a(x)$ become the classical periodic Walsh 
functions $w_n(t)$,  where $n=0, 1,2,...$ corresponds to $\alpha\in\Delta$ 
via the dyadic expansion $$n=\sum_{k=0}^\infty \alpha_{k+1}2^{k}$$

Now recall the classical Hilbert transform 
$$\Cal H f(s) = \int_{-\infty}^{\infty}{f(u)\over s-u}du $$
When 
transported to the circle coordinatized by 
$-\pi\leq \theta\leq \pi$, it becomes
$${\Cal H}f(\theta) = \int_{-\pi}^{\pi} \cot(\xi/2) f(\theta-\xi)
d\xi.
$$
  It is a complex structure on $L^2(\Bbb T)$, i.e., ${\Cal H}^2=-I$, which evidently
preserves the ordinary real form  of real-valued 
functions 
$$L^2(\Bbb T)_{\Bbb R} = \{f\in L^2(\Bbb T): \ \overline{f(t)}=f(t)\}.$$

If we set
$$\Cal H'f(s) = i\int_{-\infty}^{\infty}{f(u)
\over s+u}du,$$
then ${\Cal H'}^2=-I$ and 
$$\Cal H \Cal H' = -\Cal H' \Cal H.$$
Since any two unitary complex or 
quaternionic structures on a Hilbert space 
are mutually conjugate by a unitary 
transformation, {\it every spinor 
structure on $L^2(\Bbb T)$ can be assumed to 
start with $J_1=\Cal H$ and $J_1'=\Cal H'$}. 
Indeed, one can adjust the unitary conjugation
in such a way that 
{\it all the remaining generators are also singular integral 
operators on $\Bbb R$}
$$J_\ell f(s) =  \int_{-\infty}^{\infty}
K_\ell(s,u)f(u)du$$ 
and similarly for the $J_\ell '$. Of course, the kernels $K_\ell, K'_\ell,$ are not 
of convolution type. They are dyadic twistings of the Hilbert kernel, which, in turn,
lead to analogs of the Cauchy kernel.

\subhead 
5. Dyadic Difference Operators
\endsubhead

For any vector-valued function $f$ on $X $ define 
$$\partial_kf(x):= 
\phi_{\d^k}(x)(f(x+\d^k)-f(x))$$
where, as usual, addition is modulo 2. These 
are natural difference operators in two ways. 
Firstly, they are  natural partial derivatives 
in $X=\Bbb Z_2^\infty$  once we fix 
the motion from $0$ to $1$ (resp., from $0$ to $1$) as positive (resp., negative). Secondly, 
if we identify $X$ with the interval $(0,1)$ so that 
$x\in X$ corresponds to $t\in (0,1)$, then 
the ordinary derivative on $(0,1)$ is 
$$f'(t) =  \lim_{k\rightarrow 
\infty} 2^k\partial_k f(x).$$ 
This follows by taking 
incremental quotients of the form
$${f(t+{(-1)^{t_k}\over 2^k})-f(t)\over 
{(-1)^{t_k}\over 2^k} }=  2^k (-1)^{t_k}(f(t+{(-1)^{t_k}\over 2^k})-f(t))$$
and noting that the translation $x\mapsto x+\d^k$ in $X$, corresponds in $(0,1)$ 
to $t\mapsto t+ (-1)^{t_k} 2^{-k}$.
Equivalently,
$$
{d\over dt} = \sum_{k=0}^{\infty}\ 2^k(2\partial_{k+1}-
\partial_k).$$
This suggests a few obvious  deformations of the derivative operator, starting 
with
 $$\sum_{k=0}^{\infty}\ z^k(2\partial_{k+1}-
\partial_k),$$ 
$z\in\Bbb C$.
Another, related to 
the subject at hand, is 
obtained by  expressing
 the operators 
$\partial_k$ in terms of the the $J_k$ and $J_k'$ of the special spinor structure
$V(\mu_X,1,1)$, then replacing the operators $c_k=1$ by 
arbitrary ones. The resulting ``twisted derivative" is, for any spinor structure 
$V(\mu_X,1,\Cal C)$, 
 $${d\over d_ct} = \lim_{k\rightarrow \infty} i2^k (\phi_{\s^{k-1}}J_k+ J_kJ'_k).$$

In this 
article we will concentrate instead on 
the operators $\sum_k J_k\partial_k$ and 
$\sum_k J_k\partial_k$, or, better yet, 
$$D=\sum_{k=0}^\infty a_k\partial_k \qquad 
D'=\sum_{k=0}^\infty a_k^*\partial_k$$ 
associated to any representation $V(\mu,\nu,\Cal C)$.
We will not attempt to motivate them {\it a priori}. 
They are of course linear wherever  defined and 
anihilate constants, but their resemblance to Dirac operators does not go very far because the  $\partial_k$
do not commute with the spinor representation. 
But the following observations makes them 
worth of some attention.

 {\it For the standard Fermi-Fock 
representation the domains of $D$ and $D'$ consist of the 
constants alone} - they diverge elsewhere. However,  if $\mu_X$ is 
the Haar measure on $X\approx \Bbb T$, we have

\proclaim{Theorem 5.1} For any  
representation  $V(\mu_X,1,\Cal C)$, 
the domains of $D$ and $D'$ contain the 
algebraic span of the Walsh and the Fourier basis and, 
therefore, are dense in $L^2(X)$.
\endproclaim

In terms of the Walsh basis, the matrices of $D$ and 
$D'$ involve only $0$ and $\pm 1$ and are not 
symmetric. However they appear to be always diagonalizable. Here we shall concentrate on
what actually happens for the dyadic 
representations of \S 3, where the diagonalization can be done over $\Bbb Z$. 

Let $\g\in\Gamma$ be an infinite 
symmetric matrix of $0$'s and $1$'s such 
that all the diagonal elements 
and almost all elements in each row $\g^k$ 
are zero. The corresponding spinor  
representation on $L^2(X,\mu_X)$ is
$$\eqalign{J_kf(x) &= -i\phi_{\s^{k-1}+\g^k}(x)f(x+\d^k)\cr
J_k'f(x) &=
\phi_{\s^{k}+\g^k}(x)f(x+\d^k)\cr}$$ where $\phi_\a$ are the Walsh functions. 
The set $W_{\Bbb Z}$ of integral linear combinations 
of Walsh functions defines an integral 
structure in $L^2(X)\cong L^2(\Bbb T)$.

\proclaim{Theorem 5.2}  For any matrix $\g\in\Gamma$, the operators 
$D$ and $D'$ associated to the representation
$V(\mu_X,1,\{\phi_{\g^k}\})$ can be 
 diagonalized over $\Bbb Z$: with integral eigenvalues and 
eigenvectors in $W_{\Bbb Z}$.
\endproclaim

There is an algorithm involving only the 
matrix $\g$ to  obtain all 
the eigenvalues and   eigenvectors
of $D$ and $D'$. It goes roughly as follows.
For any positive integer $n$ let $W_n$ 
be the set of functions $X\rightarrow \Bbb C$ that depend only on the first $n$ 
components of $x$. Fix $\g\in\Gamma$ and 
define a sequence of integers $0< N_1 
< N_2 < ...$ by
$$N_k:=\max\{k, \min \{m: 
\g^1,...,\g^k\in W_{m}\}\}.$$
Then 
 $$0\subset W_{N_1}\subset W_{N_2} \subset ... $$ 
is a filtration of the space of functions on $X$ by finite-dimensional subspaces invariant under both $D$ and $D'$. 
Consider now the more 
general operators
$$D_{[k,\lambda]}:= D-\lambda\phi_{\s^{N_k}}I$$
with $k\in\Bbb N$ and $\lambda\in \Bbb C$.
Then, {\it one can 
obtain the eigenvalues and eigenvectors 
of $D_{[k,\lambda]}$ in $W_{N_k}$, recursively from the eigenvalues and
eigenvectors  of $D_{[k-1,\mu]}$ in $W_{N_{k-1}}$ for all $\mu$.}

\

In general, $D$ and $D'$ are very different operators. For example, 
the domain may be $0$ for one 
and dense for the other and there is no 
relation with the adjoints either. The 
following result essentially characterizes 
the standard spinor representations of real type.

\proclaim{Theorem 5.3} If $V(\mu_X,1,\Cal C)$ is standard 
and of real type, then 
$$ D' = TDT^{-1}$$
where
$$Tf(r) = \sqrt{\frac{d\mu(\check r)}{d\mu(r)}} {f(\check r)}.$$
\endproclaim

Finally, either operator determines the 
representation. For example, for any $V(\mu_X,1,\Cal C)$, one has

$$\eqalign{-2a_kf&=\phi_{\d^k}D(\phi_{\d^k}f)-Df\cr
-2a_k^*f&=\phi_{\d^k}D'(\phi_{\d^k}f)-D'f.\cr}
$$

\subhead 
6. Kaplansky's Division Algebras
\endsubhead

The real finite-dimensional division algebras -asociative or not- occur only in dimensions 1,2,4 
and 8. If we require them to have a multiplicative 
identity and
be {\it normed} relative to a fixed inner product,
namely, to satisfy
$$|| ab|| = ||a||\ ||b||,$$
one obtains just the usual algebras of real, complex, quaternionic and octonionic numbers. 

In [Ka], Kaplansky proved that in infinite dimensions there were no real normed division algebras, i.e. no strict analog  of the numbers above.  Of course, there are many division algebras -even asociative and commutative ones (e.g., $\Bbb R [X]$), as well as many normed algebras (since $V\otimes V \cong V$), but 
none  will satisfy both conditions simultaneously.

After noticing that weakening ``division" to, say, ``left-division", did not introduce any new algebras in finite dimensions, Kaplansky comments on his attempts to prove that the same was true in the infinite case.  But counterexamples were given in  [Cu],[R]. 

Now we can describe all such structures, that is, all bilinear operations on a real
separable Hilbert space such that $\parallel ab\parallel =\parallel a\parallel \parallel
b\parallel $ and such that for every $a\not=0$ there exists $a^{-1}$ satisfying
$a^{-1}(ab)=b$. Indeed,

\proclaim {Theorem 6.1} The left-division real normed 
algebras of countable dimension are 
parametrized  up to equivalence by the triples 
$(\mu,\nu,\Cal C)$ of Theorem 3.7.
\endproclaim

Explicitely: the product in such algebra
$A$ can be  linearly modified so as to have a left-identity $1$. If $H$ is identified with
the orthogonal complement of $1$ in $A$, so that 
$$A = H\oplus \Bbb R 1,$$
then left-multiplication on $A$ by elements of $H$ satisfies the relation 
$$(h_1h_2 +h_2h_1)a = -2<h_1,h_2>a$$
 and, therefore  $A$ becomes an orthogonal $C(H)$-module, correponding to some triple 
$(\mu,\nu,\Cal C)$ satisfying the conditions of 3.7.
Conversely, given an orthogonal
$C(H)$-module $A$ and any identification
$A\cong H\oplus{\Bbb R}$, the product
$$(h+c)\star a := ha + ca$$
where $a\mapsto ha$ is the Clifford action, satisfies the desired properties.

\

\underbar{Examples:} The Fermi-Fock representations 
yield the examples of [Cu],[R]. Letting instead $\mu_X$ be the Haar measure on $X$ and  
$$c_{2k}(x)=1, \qquad c_{4k+1}(x)=(-1)^{x_{4k+3}}, \qquad c_{4k+3}(x)=(-1)^{x_{4k+1}}$$ yields an essentially inequivalent algebra. The corresponding $C(H)$-module  is irreducible over $\Bbb C$ but splits over $\Bbb R$ as a  sum of two copies of a real irreducible module $U$. The normed algebra constructed from $U$ is therefore the simplest infinite-dimensional analog of the Octonions -if there is to be one. 

For the dyadic spinors introduced in \S 3 one can describe the resulting algebras 
purely in  dyadic terms. Let then 
$\g\in\Gamma$. 
For any  
non-negative integers $k,m,$ let
$$N_\g(k,m)=\sum_{j=0}^{k-2}(m_j +  \g^{k-1}_{j+1}
+  1)2^j  +
\sum_{j=k-1}^{\infty}(m_j +  \g^{k-1}_{j+1}
)2^j$$
$$N_\g'(k,m) =\sum_{j=0}^{k-2}(m_j +  \g^{k}_{j+1}
+  1)2^j  +
\sum_{j=k-1}^{\infty} (m_j +  \g^{k}_{j+1}
)2^j$$
where
$$k=\sum_{j\geq 0} k_j2^j, \qquad m=\sum_{j\geq 0} m_j2^j$$
($k_j,m_j\in \{0,1\}$) are the dyadic expansions of $k$ and $m$ {\it and the sums in
parenthesis are  modulo $2$}. 

\proclaim {Theorem 6.2}  Let $\g$ be an infinite matrix zeroes and ones
with finitely many
ones in each row and  such that
$$\g^k_l=\g^l_k, \qquad \g^k_k=0 $$
 for all $k,l$. 
On a real vector space $V$ with basis
$$w_0,  w_0', w_1, w_1', w_2,  w_2', ... $$
define a linear $\star_\g: A\otimes A\rightarrow A$  by 
$$\eqalign{w_0\star  w_m &= w_m\cr
 w_0\star w_m' &= w_m'\cr
w_k\star  w_m &= (-1)^{m_{k-1}}w_{N_{\g}(k,m)} \ (k\geq 1)\cr 
w_k\star  w_m' &= (-1)^{m_{k-1}}w_{N_{\g}(k,m)}' \ (k\geq 1)\cr
 w_k'\star  w_m &= -i(-1)^{m_{k-1}}w_{N_{\g}'(k,m)} 
\ (k\geq 0)\cr
 w_k'\star  w_m' &= -i(-1)^{m_{k-1}}w_{N_{\g}'(k,m)}' 
\ (k\geq 0)\cr}$$
for all $m\geq 0$. Then
 $(V,\star )$ has no zero divisors and every non-zero 
element is a left-unit. Furthermore, under the 
inner product defined by declaring $\{w_n,w_n'\}$  to be orthonormal,
$\star$ is a composition of the corresponding quadratic forms, i.e., 
$$||a\star b|| = ||a||\ ||b||$$
\endproclaim  

The algebra of the theorem corresponds to the spinor representation 
$V(\mu_X,1,\{\g^k\})$.
For example take  $\g=0$. The corresponding product $\star $ is 
antilinear in the first slot, linear 
in the second
and satisfies
$$
w_1\star  w_m = w_m$$
$$w_k\star  w_m = (-1)^{m_{k-1}} w_{m^{(k-1)}}$$
for $k\geq 2$, where $N(k,m)= m^{(k-1)} $ is the number obtained from $m$ by changing 
its first $k-1$ dyadic coefficients $m_0,..., m_{k-2}$.

\

More interesting than the products themselves 
may be their automorphism groups of various 
kinds. In particular, $\roman{Pin}(\infty)$, the Banach Lie group generated by the elements of 
unit lenght in $H$ under the Clifford product, comes with a natural unitary spin 
representation $v\mapsto J_1...J_rv$, once we fix 
one for $C(H)$. 
It satisfies
$$J_h(J_k\cdot v) = -J_{r_h(k)}\cdot J_h(v)$$
where 
$r_h:H\rightarrow H$ 
denotes the reflection with respect to the hyperplane 
$h^{\perp}$.
Therefore, $\roman{Pin}(\infty)$ acts by 
orthogonal transformations  in $H$. 
Let it act trivially on the factor $\Bbb R 1$ of $H\oplus \Bbb R 1$. If now 
$V$ is any spin representation, an identification $V = H\oplus
\Bbb R 1 $ yields two actions of $\roman{Pin}(\infty)$ on $V$,
$v\mapsto B_gv$ and $v\mapsto \Sigma_gv$, which satisfy
$$\Sigma_g(u\star v) =  B_gu \star 
\Sigma_gv. $$
In this way we obtain an inclusion 
$$1\rightarrow \roman{Pin}(\infty)\rightarrow 
\Cal G = \{(g_1,g_2)\in U(H)\times U(V): g_2(u\star v) = 
g_1u\star g_2v\}$$
for any real representation of $C(H)$. 
$\Cal G$ is reductive and comes with a unitary representation. Its specific structure 
depends very much on the equivalence class of the spinor representation but, at least in the examples treated here, their semisimple part is a classical infinite-dimensional group. Hence, by restriction, one obtains 
many irreducible unitary  representations 
of the latter. Most are ``new" and not of highest 
weight type.

For an inspiring discussion of the groups of symmetries 
asociated to the Octonions, see [B].

\enddocument
\end